\documentstyle[12pt]{article}

\date{ }
\begin{document}
\setcounter{page}{1}
\rightline{q-alg/9502018}
\rightline{CPTH-RR384.0295}
\begin{center}
\begin{Large}
The character table of the Hecke algebra $H_n(q)$ in terms of traces
of products of Murphy operators

\vspace{10pt}

$\mbox{J. Katriel}^{*}$
, B. Abdesselam and A. Chakrabarti \\
\vspace{10pt}
{\small \sl Centre de Physique Th\'eorique}\\
{\small \sl Ecole Polytechnique}\\
{\small \sl 91128 Palaiseau Cedex, France}
\end{Large}
\end{center}
\vspace*{105mm}
\ \hrulefill\ \hfill \

\begin{small}
${\ }^{*}$Permanent address: Department of Chemistry,
Technion, 32000 Haifa, Israel.
\end{small}

\newpage
\pagestyle{plain}

\vspace*{2 cm}

\begin{abstract}

The traces of the Murphy operators of the Hecke algebra $H_n(q)$,
and of products of sets of Murphy operators with non-consecutive indices,
can be evaluated by a straightforward recursive procedure.
These traces are
shown to determine all the reduced traces in this algebra, which, in turn,
determine all other traces.
To illustrate the procedure we obtain the set of reduced traces
for $H_7(q)$ - the lowest order Hecke
algebra whose character table has not hitherto been reported.
This is preceded by the presentation of an explicit algorithm for the reduction
of the trace of an arbitrary element of the Hecke algebra into a linear
combination of
traces of elements consisting of
appropriately defined disjoint cycles; and of a proof, presented in order to
make
the present article reasonably self-contained, that a
reduced trace depends only on the set of lengths of the disjoint
cycles that it consists of.

\end{abstract}

\newpage

\section{Introduction}

\hspace*{6 mm}
The present article is a sequel to \cite{Katriel}, in which the remarkable
properties
of the fundamental invariant of the Hecke algebra $H_n(q)$ were elucidated.
The sequence of sub-algebras $H_2(q)\subset H_3(q)\subset\cdots\subset H_n(q)$
gives rise to a sequence of mutually commuting fundamental invariants
$C_2,\, C_3,\, \cdots,\, C_n$. In the generic case, {\it i.e.,} when $q$ is
neither a root of unity of order up to $n$, nor zero, the eigenvalues of $C_n$,
the fundamental invariant of $H_n(q)$, fully characterize the corresponding
irreducible representations (irreps). In fact, these eigenvalues are
polynomials
in $q$ and $\frac{1}{q}$, and the coefficients of the various powers are
directly related to the structure of the Young diagram specifying the irrep.
Furthermore, the common set of eigenvectors of $C_2,\, C_3,\, \cdots ,\, C_n$,
for any {\it particular} eigenvalue of $C_n$ but all the consistent eigenvalues
of $C_2,\, C_3,\, \cdots,\, C_{n-1}$, spans the corresponding irreducible
invariant subspace. In this context it turns out that the set of
(mutually commuting) Murphy operators
$L_2=C_2,\, L_3=C_3-C_2,\, \cdots,\, L_n=C_n-C_{n-1}$ have particularly
attractive properties, the primary one being that the corresponding set of
eigenvalues characterizes in a very convenient and explicit form a specific
eigenvector within the basis specified above.

The purpose of the present article is to provide a systematic
procedure that enables the evaluation of arbitrary traces in terms
of the easily evaluated traces of products of sets of Murphy operators.
We begin by presenting an algorithm that effects the reduction of the trace of
an
arbitrary element of the Hecke algebra into a linear combination of reduced
traces. The latter consist of disjoint sequences of consecutive generators.
The fact that a reduced trace only depends on the lengths of the disjoint
sequences
of generators of which it consists was noted by King and Wybourne \cite{King}
on the basis of explicit computations for the Hecke algebras $H_2(q)$ to
$H_6(q)$,
but had in fact been shown in an unpublished doctoral dissertation by
Starkey \cite{Starkey}. To make the present article reasonably self-contained
we present in section 4 a proof of this property, formulated in three lemmas.
The direct evaluation of the traces of the reduced elements of the algebra,
without prior explicit construction of the irreps, was the main purpose of
ref. \cite{King}, which is the most comprehensive treatment of the evaluation
of
traces of elements of the Hecke algebra that we are aware of, and to which
we refer for a discussion of earlier relevant references. Further
closely relevant references are Ram \cite{ARam} and Ueno and Shibukawa
\cite{Ueno}.
However, we show in the present article that the systematic investigation of
the Murphy operators provides a
remarkably simple and convenient means for the evaluation of the reduced
traces.

\section{The Hecke algebra $H_n(q)$}

\hspace*{6 mm}
The Hecke algebra $H_n(q)$ is defined in terms of the generators
$g_1,\, g_2, \, \cdots, \, g_{n-1}$ and the relations
\begin{equation}
\label{equation:Hecke}
\begin{array}{ll}
g_i^2=(q-1)g_i+q &\;\;\; i=1,\, 2,\, \cdots,\, n-1 \\
g_ig_{i+1}g_i=g_{i+1}g_ig_{i+1} &\;\;\; i=1,\, 2,\, \cdots,\, n-2 \\
g_ig_j=g_jg_i &\;\;\; {\mbox{if }} |i-j|\geq 2
\end{array}
\end{equation}

For $q=1$ these relations reduce to the generating relations of the
symmetric group, $S_n$. In particular, $g_i$ reduces to the transposition
$(i,i+1)$.

When $q$ is neither zero nor a $k$th root of unity, $k=2,\, 3,\, \cdots, \, n$,
the irreps of $H_n(q)$ are labelled by Young diagrams with $n$ boxes
\cite{Wenzl,Vershik}.

The Murphy operators are \cite{Dipper,Murphy}
\begin{eqnarray}
\label{equation:Murphy}
L_p &=& g_{p-1}+
\frac{1}{q} g_{p-2}g_{p-1}g_{p-2}
+ \frac{1}{q^2} g_{p-3}g_{p-2}g_{p-1}g_{p-2}g_{p-3}
\\ & & {\phantom{space here}} +\cdots
+ \frac{1}{q^{p-2}} g_1g_2\cdots g_{p-2}g_{p-1}g_{p-2}\cdots g_2g_1  \nonumber
\\
&=&
\sum_{i=1}^{p-1} \frac{1}{q^{p-1-i}}g_ig_{i+1}\cdots
g_{p-2}g_{p-1}g_{p-2}\cdots
g_{i+1}g_i   \; ;\;\;\; \; \; \; \; \;
p=2,\, 3,\, \cdots,\, n \nonumber
\end{eqnarray}

Any two Murphy operators commute with one another. A state labelled by the
sequence of Young diagrams
$\Gamma_2\subset\Gamma_3\subset\cdots\subset\Gamma_n$
is an eigenstate of all the Murphy operators $L_2,\, L_3,\, \cdots,\, L_n$,
the eigenvalue of $L_i$ being \cite{Dipper,Murphy}
$$\{\Gamma_i\setminus\Gamma_{i-1}\}_q\equiv q[\kappa_i-\rho_i]_q \; ,$$
where $\kappa_i$ and $\rho_i$ are the column and row indices, respectively, of
the
box that has been added to $\Gamma_{i-1}$ to obtain $\Gamma_i$,
and $[x]_q \equiv \frac{q^x-1}{q-1}$.

The fundamental invariant of the Hecke algebra is the sum of the Murphy
operators,
$$C_n=\sum_{i=2}^n L_i \; .$$
Its eigenvalue, corresponding to an irrep $\Gamma_n$, is given by
$$\Lambda_{C_n}^{\Gamma_n} = q\sum_{(i,j)\in\Gamma_n} [j-i]_q \; ,$$
where $i$ and $j$ are the row and column indices, respectively, of the boxes
comprising the Young diagram that corresponds to $\Gamma_n$.

It will be convenient to introduce the shorthand notation
$$ (g_i)^{(\ell)}\equiv g_ig_{i+1}\cdots g_{i+\ell-1}\cdots g_{i+1}g_i$$
and
$$(g_i)_{\ell}\equiv g_ig_{i+1}\cdots g_{i+\ell-1} \; .$$
Using this notation
$$L_p=\sum_{i=1}^{p-1} \frac{1}{q^{p-1-i}} (g_i)^{(p-i)}=
\sum_{j=1}^{p-1} \frac{1}{q^{j-1}} (g_{p-j})^{(j)} \; .$$

The  generalized braiding relation
$$g_ig_{i+1}\cdots g_{j-1}g_jg_{j-1}\cdots g_{i+1}g_i=
g_jg_{j-1}\cdots g_{i+1}g_i g_{i+1}\cdots g_{j-1}g_j$$
is easily proved and sometimes useful.

Some remarkable and handy identities that are easily derived are
$$g_i(g_j)_{m-j+1}=(g_j)_{m-j+1}g_{i-1}\;\;\;\;\;\;{\mbox{for }}\; i\leq j-2\;
;
\; j+1\leq i\leq m \; ;\; i\geq m+2$$
$$g_i(g_j)^{(m-j+1)}=(g_j)^{(m-j+1)}g_i\;\;\;\;\;\;{\mbox{for }}\; i\leq j-2\;
;
\; j+1\leq i\leq m-1\; ; \; i\geq m+2$$
\begin{equation}
\label{equation:Lpp}
L_{p+1}=
  \frac{1}{(q-1)q^{p-1}}\Big((g_pg_{p-1}\cdots g_1)(g_1\cdots
g_{p-1}g_p)-q^p\Big)
\end{equation}
 From the last identity it follows that
$$L_{p+1}=\frac{1}{q}g_pL_pg_p + g_p\; .$$

\section{Reduction of the trace of an arbitrary element of the Hecke algebra}

\hspace*{6 mm}
We shall refer to a product of $\ell$ consecutive generators of the Hecke
algebra,
$(g_i)_{\ell}$,
as a connected sequence of length $\ell$, or sometimes, generalizing the
corresponding
notion from the context of the symmetric group, as a cycle of length $\ell+1$.
Two connected sequences are disjoint if the sequence with a lower bottom
element has a
top element whose index
is less than and not
adjacent to that of the bottom element of the other sequence.
Obviously, two disjoint sequences commute with one another.
An element of the Hecke algebra consisting of disjoint cycles will be referred
to as a reduced element.

King and Wybourne \cite{King} have demonstrated the reduction of the traces
of arbitrary elements in the Hecke algebras $H_2(q)$ to $H_6(q)$ into linear
combinations of reduced traces, {\it i.e.,} traces of reduced elements.
In the present section we provide an algorithm for the reduction of the trace
of
an arbitrary element of $H_n(q)$.

The following observation will be needed. A trace in which each generator
appears at most
once can be transformed into the trace of the same set of generators, ordered
by
increasing indices. In the latter form the disjoint connected
sequences are readily identified.
At the risk of stating the obvious
we note that the transformation we refer to can be carried out recursively, the
$\ell$th step being
$$tr\Big((\phi)_{\ell-1}Dg_{\ell}E\Big)=
tr\Big(D(\phi)_{\ell-1}g_{\ell}E\Big)=
tr\Big((\phi)_{\ell-1}g_{\ell}ED\Big)$$
where $D$ and $E$ contain no generator with index less then $\ell+1$, and all
higher
generators at most once, and
where $(\phi)_{\ell-1}$ is an ordered product of the generators $g_1$, $g_2$,
$\cdots$,
$g_{\ell-1}$, containing each one at most once.

An algorithm that achieves the reduction of the
trace of an arbitrary product of generators
into a linear combination of reduced traces,
in a finite number of
steps, is formulated as follows.
Assume that we have already transformed the trace of the original product of
generators
into a linear combination of traces of the form
\begin{equation}
\label{equation:arbi}
tr\Big((\phi)_{\ell-1}g_{\ell}Ag_{\ell}B\Big)
\end{equation}
where $A$ is an arbitrary product of arbitrary numbers of generators with
indices not
less then $\ell+1$, and
$B$ is an arbitrary product of arbitrary numbers of generators with indices not
less
then $\ell$.

\begin{itemize}

\item
If $g_{\ell+1}$ does not appear within $A$, the two $g_{\ell}$ factors can be
brought
together, and reduced.

\item
If $g_{\ell+1}$ appears within $A$ precisely once, the two $g_{\ell}$ factors
can be carried
to its sides and the braiding relation used to reduce the number of $g_{\ell}$
factors
by one.

\item
If $g_{\ell+1}$ appears within $A$ more than once, the leftmost segment of the
form $g_{\ell+1}A^{\prime}g_{\ell+1}$ within $A$ should be treated following
the
analogue of the appropriate one of the steps that are presently being
described.
Once the number of $g_{\ell+1}$ factors between the first two $g_{\ell}$
factors
has been reduced to zero or one, the appropriate mode of reducing the number
of $g_{\ell}$ factors should be applied.

\end{itemize}

Each one of the traces thus obtained can be brought
back to the form of eq. \ref{equation:arbi} with redefined $A$ and $B$,
possibly by appropriate cyclic
transformations within each trace, or transposition of commuting factors.

The whole procedure should continue until all the traces generated contain
$g_{\ell}$
at most once, at which stage each trace should be brought to the form
$$tr\Big((\phi)_{\ell}g_{\ell+1}Ag_{\ell+1}B\Big) \; ,$$
with appropriately redefined $A$ and $B$, now containing factors with indices
at least $\ell+2$ and $\ell+1$, respectively.

Note that the starting point of the iteration is a trace of the form
$tr\Big(g_1Ag_1B\Big)$ where we have interpreted $(\phi)_0$ to be an empty
sequence, that is
equal to the identity, and where $A$ does not contain $g_1$.

Since the number of generators is finite the algorithm terminates in a finite
number of iterations.
The resulting expression is a linear combination of reduced traces.

\section{Traces of reduced elements with common cycle structure}

\hspace*{6 mm}
We shall proceed by proving the following three {\it Lemmas}.

\noindent
{\sc Lemma 1:}  {\it The trace of a connected sequence of a given length, in
any
particular irrep, depends only on the length of the sequence.}

\noindent
{\sc proof:} First, note that for obvious reasons the first index in a
connected
sequence of length $\ell$ must satisfy $1\leq i\leq n-\ell$.

\noindent
Now,

\begin{eqnarray}
& & tr\Big(g_{i+\ell} (g_i)_{\ell} g_{i+\ell}\Big)=
(q-1)tr\Big((g_i)_{\ell+1}\Big) +q\, tr\Big((g_i)_{\ell}\Big)
\nonumber \\
& &{\phantom{space}}
=tr(g_ig_{i+1}\cdots g_{i+\ell}g_{i+\ell-1}g_{i+\ell})=
tr(g_ig_{i+1}\cdots g_{i+\ell-1}g_{i+\ell}g_{i+\ell-1})= \nonumber\\
& &{\phantom{space}}
=tr(g_{i+\ell-1}g_ig_{i+1}\cdots g_{i+\ell-1}g_{i+\ell})=
tr(g_ig_{i+1}\cdots g_{i+\ell-1}g_{i+\ell-2}g_{i+\ell-1}g_{i+\ell})=
\nonumber \\ & &{\phantom{space}}
=tr(g_ig_{i+1}\cdots g_{i+\ell-2}g_{i+\ell-1}g_{i+\ell-2}g_{i+\ell})=
tr(g_ig_{i+1}\cdots
g_{i+\ell-2}g_{i+\ell-3}g_{i+\ell-2}g_{i+\ell-1}g_{i+\ell})
\nonumber \\ & & {\phantom{space}}      =\; \cdots\; =
\nonumber \\ & &{\phantom{space}}
=tr(g_{i+1}g_ig_{i+1}g_{i+2}\cdots g_{i+\ell-1}g_{i+\ell})=
tr(g_ig_{i+1}g_ig_{i+2}\cdots g_{i+\ell-1}g_{i+\ell})=
\nonumber \\ & &{\phantom{space}}
=(q-1)tr\Big((g_i)_{\ell+1}\Big)
+ q\, tr\Big((g_{i+1})_{\ell}\Big)
\;\;\;\;\;\;\; ;\;\;\; i=1,\, 2,\, \cdots,\, n-\ell-1\,.
\nonumber
\end{eqnarray}

\noindent
Hence,
$$tr\Big((g_i)_{\ell}\Big)
= tr\Big((g_{i+1})_{\ell}\Big)
\;\;\;\;\;\;\; ;\;\;\; i=1,\, 2,\, \cdots,\, n-\ell-1\, ,$$
which is a statement of the {\it Lemma}.

\noindent
{\sc Lemma 2:}  {\it The trace of a product of two disjoint connected
sequences,
in any particular irrep, depends only on the lengths of the two sequences.}

\noindent
{\sc proof:} Consider the trace
$$tr\Big((g_i)_{\ell} \; (g_j)_m\Big)_{\Gamma} \; ,$$
where $i+\ell<j$.
First, note that using the argument of {\it Lemma 1} the indices of each
connected
sequence in this trace can be shifted as follows: $i$ can be shifted to any
value
between 1 and $j-\ell-1$ without changing $j$. Independently, $j$ can be
shifted to
any value between $i+\ell+1$ and $n-m-1$ without changing $i$.

\noindent
Now, write the above trace in the equivalent form
$$tr\Big((g_1)_{\ell} \; (g_{n-m})_{m}\Big)_{\Gamma} \; ,$$
{\it i.e.,} shift $i$ to the bottom and $j$ to the top.
We now show that the trace remains unchanged if the two connected sequences are
transposed, so that the sequence of length $\ell$ would consist of generators
with higher indices than the sequence of length $m$.
Note the automorphism $g_i\leftrightarrow {\tilde{g}}_i\equiv g_{n-i}$,
$i=1,\, 2,\, \cdots,\, n-1$, and consider some irrep $\Gamma:\; g_i\mapsto
D_i$.
Obviously, ${\tilde{\Gamma}} : \; g_i\mapsto {\tilde D}_i\equiv D_{n-i}$
is an equivalent irrep.
Let $U$ be the
transformation matrix from the first to the second irrep, {\it i.e.,}
$D_i=U {\tilde D}_i U^{-1}=U D_{n-i}U^{-1}$.
It follows that
\begin{eqnarray}
& & tr\Big((g_1)_{\ell} \; (g_{n-m})_{m}\Big)_{\Gamma}=
tr(D_1D_2\cdots D_{\ell}\; D_{n-m}D_{n-m+1}\cdots D_{n-1}) =
\nonumber \\ & &
=tr(UD_{n-1}U^{-1}UD_{n-2}U^{-1}\cdots UD_{n-\ell}U^{-1}\;
UD_{m}U^{-1} UD_{m-1}U^{-1} \cdots UD_1U^{-1})= \nonumber \\ & &
=tr(UD_{n-1}D_{n-2}\cdots D_{n-\ell}\;
D_{m}D_{m-1} \cdots D_1U^{-1})= \nonumber \\ & &
=tr(g_{n-1}g_{n-2}\cdots g_{n-\ell} \; g_{m}g_{m-1}\cdots g_{1})_{\Gamma}
=tr\Big((g_1)_m \; (g_{n-\ell})_{\ell}\Big)_{\Gamma}
\nonumber
\end{eqnarray}
which proves the {\it Lemma}.

\noindent
{\sc Lemma 3:}  {\it The trace of a product of any number of disjoint connected
sequences,
in any particular irrep, depends only on the set of lengths of the different
connected sequences.}

\noindent
{\sc proof:} We shall proceed by an induction on the number $\kappa$ of
disjoint connected
sequences. From {\it Lemma 2} it follows that the present {\it Lemma} is true
for
$\kappa=2$. Assume that it is true for $\kappa-1$.

We shall denote by
$\langle\cdots\rangle_{\Gamma_i\subset\Gamma_{i+1}\subset\cdots\subset\Gamma_j}$
the diagonal matrix element of $\cdots$ corresponding to any one of the basis
vectors specified by
$\Gamma_i\subset\Gamma_{i+1}\subset\cdots\subset\Gamma_j$.
The suppression of the subsequence
$\Gamma_2\subset\Gamma_3\subset\cdots\subset\Gamma_{i-1}$
in the symbol specifying the basis vector implies that the matrix element in
question does not depend on that subsequence.
Summation over
$(\Gamma_i\subset\Gamma_{i+1}\subset\cdots\subset\Gamma_{j-1})\subset\Gamma_j$
means that
$\Gamma_i\subset\Gamma_{i+1}\subset\cdots\subset\Gamma_{j-1}$
obtain all the combinations of values consistent with the given $\Gamma_j$.

\noindent
Given
$$\tau\equiv tr\left((g_{i_1})_{\ell_1} \;
(g_{i_2})_{\ell_2}\; \cdots
(g_{i_{\kappa-1}})_{\ell_{\kappa-1}} \;
(g_{i_{\kappa}})_{\ell_{\kappa}}\right)_{\Gamma_n}$$
where $1\leq i_1$, $i_1+\ell_1<i_2$, $i_2+\ell_2<i_3$, $\cdots$,
$i_{\kappa-1}+\ell_{\kappa-1}<i_{\kappa}$, $i_{\kappa}+\ell_{\kappa}\leq n$,
we begin by shifting the last connected sequence to the top, obtaining
\begin{eqnarray}
\tau & =& tr\left((g_{i_1})_{\ell_1} \;
(g_{i_2})_{\ell_2} \; \cdots
(g_{i_{\kappa-1}})_{\ell_{\kappa-1}} \;
(g_{n-\ell_{\kappa}-1})_{\ell_{\kappa}} \right)_{\Gamma_n}= \nonumber \\
& & =
\sum_{(\Gamma_2\subset\Gamma_3\subset\cdots\subset\Gamma_{n-1})\subset\Gamma_n}
\left\langle (g_{i_1})_{\ell_1} \;
(g_{i_2})_{\ell_2} \; \cdots
(g_{i_{\kappa-1}})_{\ell_{\kappa_1}} \;
(g_{n-\ell_{\kappa}-1})_{\ell_{\kappa}}
\right\rangle_{\Gamma_2\subset\Gamma_3\subset\cdots\subset\Gamma_{n-1}\subset\Gamma_n}
\nonumber
\end{eqnarray}

The operation of $g_i$ on the state labelled by
$$\Gamma_2\subset\Gamma_3\subset\cdots\subset\Gamma_{i-1}
\subset\Gamma_i\subset\Gamma_{i+1}
\subset\cdots\subset\Gamma_{n}$$
can only affect the Young diagram $\Gamma_{i-1}$,
\cite{Pan}, so that \hfill\break
$\langle(g_i)_{\ell}\rangle_{\Gamma_2\subset\Gamma_3\subset\cdots\subset\Gamma_{n}}=
\langle(g_i)_{\ell}\rangle_{\Gamma_{i-1}\subset\Gamma_{i}\subset\cdots\subset
\Gamma_{i+\ell-2}}$.
Therefore,
\begin{eqnarray}
& & \tau=
\sum_{(\Gamma_{n-\ell_{\kappa}-2}\subset\Gamma_{n-\ell_{\kappa}-1}\subset\cdots\subset
\Gamma_{n-1})\subset\Gamma_n} \nonumber \\ & &
\left\{ \sum_{(\Gamma_2\subset\Gamma_3\subset\cdots\subset
\Gamma_{n-\ell_{\kappa}-3} )\subset\Gamma_{n-\ell_{\kappa}-2}}
\left\langle(g_{i_1})_{\ell_1} \;
(g_{i_2})_{\ell_2} \; \cdots
(g_{i_{\kappa-1}})_{\ell_{\kappa-1}}
\right\rangle_{\Gamma_2\subset\Gamma_3\subset\cdots
\subset\Gamma_{n-\ell_{\kappa}-3}\subset\Gamma_{n-\ell_{\kappa}-2}}\right\}
\nonumber \\ & &
{\phantom{empty line}} \nonumber \\ & &
{\phantom{here we put space to make it look much nicer}}
\times\left\langle
(g_{n-\ell_{\kappa}-1})_{\ell_{\kappa}}
\right\rangle_{\Gamma_{n-\ell_{\kappa}-2}\subset\Gamma_{n-\ell_{\kappa}-1}\subset
\cdots\subset
\Gamma_{n-1}\subset\Gamma_n}
\nonumber \\ & &
=\sum_{(\Gamma_{n-\ell_{\kappa}-2}\subset\Gamma_{n-\ell_{\kappa}-1}\subset\cdots\subset
\Gamma_{n-1})\subset\Gamma_n}
\nonumber \\ & & {\phantom{beauty}}
tr\left(
(g_{i_1})_{\ell_1} \;
(g_{i_2})_{\ell_2} \; \cdots
(g_{i_{\kappa-1}})_{\ell_{\kappa-1}}
\right)_{\Gamma_{n-\ell_{\kappa}-2}}
\left\langle
(g_{n-\ell_{\kappa}-1})_{\ell_{\kappa}}
\right\rangle_{\Gamma_{n-\ell_{\kappa}-2}\subset\Gamma_{n-\ell_{\kappa}-1}\subset
\cdots\subset
\Gamma_{n-1}\subset\Gamma_n} \nonumber
\end{eqnarray}
Now, within the trace of the product of the first $\kappa-1$ disjoint
connected sequences
the induction hypothesis allows arbitrary permutations of
disjoint sequences as well as arbitrary
shifts of indices that respect the disjointness conditions. Performing such a
combination of permutations and shifts we establish the {\it Lemma}, except for
permutations involving the top connected sequence.
Permutations of that kind can be performed using the argument of the proof of
{\it Lemma 2} to transpose the top and bottom connected sequences in the
original form of $\tau$.
This proves the {\it Lemma}.

The simplest non-trivial applications of {\it Lemmas} 1, 2, and 3 are the
identities
$tr(g_1)=tr(g_2)$, $tr(g_1g_2 \; g_4)=tr(g_1 \; g_3g_4)$, and
$tr(g_1 \; g_3g_4 \; g_6)=tr(g_1\; g_3 \; g_5g_6)$, respectively.

\section{Some basic reduction formulas}

\hspace*{6 mm}
In the following section we demonstrate that the derivation of expressions for
arbitrary reduced traces in terms of traces of products of Murphy operators
requires the reduction of traces of the form
$$tr\Big((g_1)_{\ell_1} (g_{m_2})_{\ell_2}\cdots
(g_{m_{\kappa}})_{\ell_{\kappa}}
(g_i)^{(p)}\Big) \; ,$$
where $m_j=\ell_1+\ell_2+\cdots+\ell_{j-1}+ j$ and where
$(g_1)_{\ell_1} (g_{m_2})_{\ell_2}\cdots (g_{m_{\kappa}})_{\ell_{\kappa}}$
are $\kappa$ disjoint cycles. The index $i$ can assume any value,
but $i+p-1>m_{\kappa}+\ell_{\kappa}-1$. In fact, the complete range of cases is
covered by $1\leq i\leq \ell_1+2$, $\kappa=1,\, 2,\, \cdots$. This follows
from the fact that for $i\geq m_j$ the first $j-1$ cycles commute with all
the generators comprising $(g_i)^{(p)}$ and, consequently, remain passive
``spectators''. This feature will play a central role in the following section.
The special cases of this
expression that will be required in section 6, where we consider all the
reduced
traces appearing in $H_2(q)$ up to $H_7(q)$, are dealt with in the present
section.
The general procedure is presented in the Appendix.

The identity
\begin{equation}
\label{equation:tr1}
tr\Big((g_1)^{(p)}\Big)= q^{p-1}
\sum_{i=0}^{p-1} {{p-1}\choose i} \left(\frac{q-1}{q}\right)^i \;
tr\Big((g_1)_{i+1}\Big)
\end{equation}
is obtained by straightforward application of the defining relations of the
Hecke
algebra, eq. \ref{equation:Hecke}, and of the properties of the trace, in
particular
$tr(AB)=tr(BA)$.

\vfill\eject

The recurrence relation
\begin{eqnarray}
\label{equation:l3l}
& & tr\Big((g_1)_p(g_1)^{(\ell)}\Big)=
f_{2p+1}tr\Big((g_1)_p(g_{p+1})^{(\ell-p)}\Big) \nonumber \\ & & {\phantom{some
space}}
+(q-1)\sum_{\ell=1}^p q^{\ell}f_{2(p-\ell)+1}
tr\Big((g_1)_{\ell-1}\; (g_{\ell+1})_{p-\ell}(g_{p+1})^{(\ell-p)}\Big)
\end{eqnarray}
where $f_p=\frac{q^p-(-1)^p}{q+1}$ and $\ell>p$, is demonstrated in the
Appendix.

The simplest special case is
\begin{eqnarray*}
tr\Big(g_1(g_1)^{(\ell)}\Big)&=&(q-1)tr\Big((g_1)^{(\ell)}\Big)
+q\, tr\Big(g_1(g_2)^{(\ell-1)}\Big) \\
&=&\frac{q^3+1}{q+1}tr\Big(g_1(g_2)^{(\ell-1)}\Big)+
q(q-1) tr\Big((g_2)^{(\ell-1)}\Big)
\end{eqnarray*}
and the next is
\begin{eqnarray}
\label{equation:L2l}
& & tr\Big((g_1)_2(g_1)^{(\ell)}\Big)=
\frac{q^5+1}{q+1}tr\Big((g_1)_2(g_3)^{(\ell-2)}\Big) \nonumber \\ & &
+q(q-1)\frac{q^3+1}{q+1}\,
tr\Big(g_1(g_2)^{(\ell-2)}\big)+q^2(q-1)tr\Big(g_1(g_3)^{(\ell-2)}\Big)
\end{eqnarray}

The recurrence relation
$$tr\Big((g_1)_{\ell-1}(g_{\ell})^{(k-\ell+1)}\Big)=
(q-1)\, tr\Big((g_1)_{\ell}(g_{\ell+1})^{(k-\ell)}\Big)+
q\, tr\Big((g_1)_{\ell-1}(g_{\ell})^{(k-\ell)}\Big)$$
is straightforwardly derived. It can be used to obtain the explicit expression
\begin{equation}
\label{equation:l1l}
tr\Big((g_1)_{\ell-1}(g_{\ell})^{(m+1)}\Big)=
q^m\sum_{i=0}^m {m\choose
i}\left(\frac{q-1}{q}\right)^itr\Big((g_1)_{\ell+i}\Big)\; .
\end{equation}

Furthermore,
\begin{equation}
\label{equation:ll}
tr\Big((g_1)_{\ell}(g_{\ell})^{(m+1)}\Big)=
q^m\sum_{i=0}^m\left\{(q-1){m\choose
i}+\frac{q}{q-1}{{m-1}\choose{i-1}}\right\}
\left(\frac{q-1}{q}\right)^i tr\Big((g_1)_{\ell+i}\Big)
\end{equation}

Finally, when $k>2$
\begin{eqnarray}
\label{equation:g1gi}
&&tr\Big((g_1)_3(g_2)^{(k)}\Big)=\frac{q^5+1}{q+1}tr\Big((g_1)_3(g_4)^{(k-2)}\Big)
+q(q-1)\frac{q^3+1}{q+1}tr\Big((g_1)_2(g_3)^{(k-2)}\Big) \nonumber
\\&&{\phantom{here we really need a lot of space}}
+q^2(q-1)tr\Big(g_1\, g_3(g_4)^{(k-2)}\Big)
\end{eqnarray}
Note that the expressions reduced in eqs. \ref{equation:L2l} and
\ref{equation:g1gi} involve the same number of overlapping generators, and,
consequently, have the same structure. This property is fully elucidated in the
Appendix.

\section{Reduced traces in terms of traces of products of Murphy operators}

\hspace*{6 mm}
It is a simple matter to construct
the (diagonal) representation matrices of the Murphy operators, from which
their
traces, as well as those of products of arbitrary subsets of Murphy operators,
can readily be obtained.
In fact, we show in section 7 that the traces of products of Murphy operators
can
be directly evaluated by a very simple recursive procedure, so that the actual
representation matrices need not be constructed.

In the present section we discuss the evaluation of reduced traces in terms of
traces of products of Murphy operators. For the sake of clarity of the
presentation
we proceed by considering cases of ascending complexity, thereby introducing
the procedures required in their simplest possible form. Our
guideline in the present section is to develop explicit expressions for all
the reduced traces appearing in Hecke algebras up to $H_7(q)$, inclusive.
Once this is achieved the generalization to arbitray Hecke algebras
becomes obvious, in principle.
Note, however, that the expression obtained for any particular reduced trace
in terms of appropriate traces of Murphy operator products, within the lowest
Hecke algebra within which that reduced trace appears,
remains valid for the traces
of that particular type for all $H_n(q)$.

Using eqs. \ref{equation:Murphy} and \ref{equation:tr1} and the binomial
identity $\sum_{i=1}^{p-1-j} {{p-1-i}\choose j}={{p-1}\choose{j+1}}$ we obtain
$$tr(L_k)=\sum_{i=0}^{k-2} {{k-1}\choose{i+1}}
\left(\frac{q-1}{q}\right)^{i}tr\Big((g_1)_{i+1}\Big)$$
Inverting, we obtain explicit expressions for all the reduced traces consisting
of single cycles (singly connected sequences) in terms of traces of Murphy
operators
\begin{equation}
\label{equation:simply}
tr\Big((g_1)_{k-1}\Big)=\left(\frac{q}{q-1}\right)^{k-2}
\sum_{i=0}^{k-2} (-1)^i{{k-1}\choose i} tr(L_{k-i})
\end{equation}

To proceed, we formulate the following recursive procedure for the evaluation
of reduced traces with $\kappa+1$ disjoint cycles, starting from the
expressions
for appropriate reduced traces  with $\kappa$ disjoint cycles.
Let $\Psi_{(\ell_1,\ell_2,\cdots,\ell_{\kappa})}$
be a product of $\kappa$ disjoint cycles of lengths
$\ell_1,\ell_2,\cdots,\ell_{\kappa}$. For convenience we assume that
$\ell_1\leq \ell_2\leq \cdots\leq \ell_{\kappa}$ and that the representative
product
$\Psi_{(\ell_1,\ell_2,\cdots,\ell_{\kappa})}$ is chosen with minimal spacings
among
consecutive cycles, {\it i.e.,}
$\Psi_{(\ell_1,\ell_2,\cdots,\ell_{\kappa})}=
(g_1)_{\ell_1}(g_{\ell_1+2})_{\ell_2}(g_{\ell_1+\ell_2+3})_{\ell_3}\cdots
(g_{\ell_1+\ell_2+\cdots+\ell_{\kappa-1}+\kappa})_{\ell_{\kappa}}$.
Let us assume that
$\Psi_{(\ell_1,\ell_2,\cdots,\ell_{\kappa})}$ has already been expressed as
a linear combination of the traces of a certain set of Murphy operators
\begin{equation}
\label{equation:psi}
tr(\Psi_{(\ell_1,\ell_2,\cdots,\ell_{\kappa})})=\sum_{(i_1,i_2,\cdots,i_{\kappa})}
C_{i_1,i_2,\cdots,i_{\kappa}}^{(\ell_1,\ell_2,\cdots,\ell_{\kappa})}
tr(L_{i_1}L_{i_2}\cdots L_{i_{\kappa}})
\end{equation}
where $0\leq i_1\leq i_2\leq \cdots\leq i_{\kappa}$ and where $L_0\equiv 1$.
To obtain the reduced traces with $\kappa+1$ cycles of lengths
$\ell_1,\ell_2,\cdots,\ell_{\kappa},\ell_{\kappa+1}$, where, without loss of
generality,
we assume that $\ell_{\kappa+1}\geq\ell_{\kappa}$, we note

\noindent
{\sc Lemma 4:}
$tr(\Psi_{(\ell_1,\ell_2,\cdots,\ell_{\kappa})}L_m)=\sum_{(i_1,i_2,\cdots,i_{\kappa})}
C_{i_1,i_2,\cdots,i_{\kappa}}^{(\ell_1,\ell_2,\cdots,\ell_{\kappa})}
tr(L_{i_1}L_{i_2}\cdots L_{i_{\kappa}}L_m)$,
where $m>\ell_1+\ell_2+\cdots+\ell_{\kappa}+\kappa$.

\noindent
{\sc proof:} Since $L_m$ commutes with all the generators comprising
$\Psi_{(\ell_1,\ell_2,\cdots,\ell_{\kappa})}$, all the manipulations leading to
the
identity \ref{equation:psi} can be carried out through $L_m$.

The remaining task is the reduction of
$tr(\Psi_{(\ell_1,\ell_2,\cdots,\ell_{\kappa})}L_m)$
into a linear combination of reduced traces. This is effected using the
recurrence
relations presented in section 5 and in the Appendix.
Assuming that all traces with up to $\kappa$ disjoint cycles, as well as all
traces with $\kappa+1$ disjoint cycles such that the first $\kappa$ cycles
preceed
$\Psi_{(\ell_1,\ell_2,\cdots,\ell_{\kappa})}$
in a lexicographic ordering, have been determined, we begin with
$m=\ell_1+\ell_2+\cdots+\ell_{\kappa}+\kappa+\ell_{\kappa}+1$ to obtain
the new trace with $\ell_{\kappa+1}=\ell_{\kappa}$, {\it i.e.,}
$tr(\Psi_{(\ell_1,\ell_2,\cdots,\ell_{\kappa},\ell_{\kappa})})$.
We then increase $m$ to obtain the reduced traces with higher
$\ell_{\kappa+1}$.

Thus, from eq. \ref{equation:simply} we obtain
\begin{equation}
\label{equation:LkLm}
tr\Big((g_1)_{k-1} L_m\Big)=\left(\frac{q}{q-1}\right)^{k-2} \sum_{i=0}^{k-2}
(-1)^i
{{k-1}\choose{i}}
tr(L_{k-i} \; L_m) \; .
\end{equation}
where $m>k+1$. Note that all the terms on the right hand side are traces
of products of Murphy operators.

To proceed, the left hand side of eq. \ref{equation:LkLm} has to be expressed
as a
linear combination of reduced traces. By a judicious choice of $m$ precisely
one of these will be of a new type. That trace would be the one that the
resulting
identity would allow evaluating.

The simplest case is $tr(g_1L_m)$, for which eq. \ref{equation:LkLm} is
trivial,
yielding $tr(g_1L_m)=tr(L_2L_m)$.
To evaluate that trace we need
the three identities
$$tr\Big(g_1 \; (g_3)^{(p-2)}\Big)= q^{p-3}
\sum_{i=0}^{p-3} {{p-3}\choose i} \left(\frac{q-1}{q}\right)^i \;
tr\Big(g_1\; (g_3)_{i+1}\Big) \; ,$$
which, like eq. \ref{equation:LkLm}, follows from the fact that $g_1$ commutes
with all
the generators comprising $(g_3)^{(p-2)}$,
$$tr\Big(g_1 (g_2)^{(p-1)}\Big)= q^{p-2}
\sum_{i=0}^{p-2} {{p-2}\choose i} \left(\frac{q-1}{q}\right)^i\;
tr\Big((g_1)_{i+2}\Big)
$$
which is a special case of eq. \ref{equation:l1l}, and
$$tr\Big(g_1(g_1)^{(p)}\Big)= q^{p-1}
\sum_{i=0}^{p-1} \left\{ (q-1){{p-1}\choose i}
+\frac{q}{q-1}{{p-1}\choose{i-1}}\right\}
 \left(\frac{q-1}{q}\right)^i\; tr\Big((g_1)_{i+1}\Big)$$
which is a special case of eq. \ref{equation:ll}.

\noindent
Thus, for $m\geq 4$
\begin{eqnarray}
\label{equation:L2Lk}
tr(g_1 L_m) &=& \sum_{j=0}^{m-4}\left(\frac{q-1}{q}\right)^j
{{m-3}\choose {j+1}} tr\Big(g_1 \; (g_3)_{j+1}\Big) \\ &+&
\sum_{j=0}^{m-2} \left\{ \frac{2q}{q-1} {{m-3}\choose{j-1}}+
(q-1){{m-2}\choose{j}}\right\}\left(\frac{q-1}{q}\right)^j
tr\Big((g_1)_{j+1}\Big) \nonumber
\end{eqnarray}
In fact, for $m=3$ eq. \ref{equation:L2Lk} reduces to
$$tr(g_1 L_3)=(q-1)tr(g_1)+\left(q+\frac{1}{q}\right)tr(g_1g_2) \; ,$$
that provides no new information. Expressing the traces on the right-hand-side
by means of $tr(L_2)$ and $tr(L_3)$ we obtain the identity
$$(q-1)tr(L_2L_3)=(q^2+1)tr(L_3)-(q+1)^2tr(L_2) \; .$$

Eq. \ref{equation:L2Lk} enables the systematic evaluation of all reduced
traces of the form\hfill\break
$tr\Big(g_1\, (g_3)_j\Big) \;\; ; \; \; j=1,\, 2,\, \cdots $,
{\it i.e.,} the traces of all elements consisting of a cycle of length two and
a disjoint cycle of arbitrary length.

As an example consider
\begin{equation}
\label{equation:g1L4}
tr(g_1L_4)=tr(g_1\,
g_3)+\frac{q-1}{q}(q+\frac{1}{q})tr\Big((g_1)_3\Big)+
2(q-1+\frac{1}{q})tr\Big((g_1)_2\Big)+(q-1)tr(g_1)
\end{equation}
which enables the evaluation of the non-simply-connected element $tr(g_1\,g_3)$
consisting of two disjoint cycles of unit length, {\it i.e.,}
\begin{equation}
\label{equation:g1g3}
tr(g_1g_3)=tr(L_2L_4)-\frac{q^2+1}{q-1}tr(L_4)+
\frac{(q+1)^2}{q-1}tr(L_3)-\frac{2q}{q-1}tr(L_2)\; .
\end{equation}

As a further example consider $tr(g_1 L_5)$, that can be expressed
in terms of $tr(g_1\, g_3)$, $tr(g_1\, g_3g_4)$, and traces of simply connected
terms.
Thus,
\begin{eqnarray}
tr(g_1L_5)&=&2\, tr(g_1\, g_3)+\left(\frac{q-1}{q}\right)tr(g_1\, g_3g_4)
\nonumber \\
&+&\left(\frac{q-1}{q}\right)^2\left(q+\frac{1}{q}\right)tr(g_1g_2g_3g_4)
+\left(\frac{q-1}{q}\right)\left(3q-2+\frac{3}{q}\right)tr(g_1g_2g_3)
\nonumber\\
&+&\left(3q-4+\frac{3}{q}\right)tr(g_1g_2)+(q-1)tr(g_1) \nonumber
\end{eqnarray}
Since $tr(g_1\, g_3)$, as well as all the simply-connected traces, have already
been
evaluated, $tr(g_1 L_5)$ provides the next non-simply connected term,
$tr(g_1\, g_3g_4)$.

Next, we consider $tr(g_1g_2 L_m)$, for which eq. \ref{equation:LkLm} yields
$$tr(g_1g_2 L_m)=\frac{q}{q-1}\Big(tr(L_3L_m)-2tr(L_2L_m)\Big) \; .$$
Since the right hand side consists of traces of products of Murphy operators,
the remaining task is the reduction of the left hand side.
We obtain
\begin{eqnarray*}
tr(g_1g_2L_m)&=&\frac{1}{q^{m-2}}tr\Big((g_1)_2(g_1)^{(m-1)}\Big)
+\frac{1}{q^{m-3}}tr\Big((g_1)_2(g_2)^{(m-2)}\Big) \\
&+& \frac{1}{q^{m-4}}tr\Big((g_1)_2(g_3)^{(m-3)}\Big)
+\sum_{i=4}^{m}\frac{1}{q^{m-1-i}}tr\Big((g_1)_2(g_4)^{(m-i)}\Big)
\end{eqnarray*}
where the first,
second and third terms on the right hand side are special cases of
eqs. \ref{equation:l3l}, \ref{equation:ll}
and \ref{equation:l1l}, respectively.
All these terms reduce into single cycle traces, with the exception of the
contribution
$tr\Big(g_1(g_3)^{(m-3)}\Big)$ to the first term. Note, however, that this
trace was
already evaluated in terms of $tr(L_2L_m)$.
Each summand in the fourth term can be reduced into two-cycle traces. Thus,
using eq. \ref{equation:tr1} and the fact that $(g_1)_2$ commutes with all the
generators comprising $(g_4)^{(j)}$,
$$tr\Big((g_1)_2 (g_4)^{(j)}\Big)=
q^{j-1}\sum_{i-0}^{j-1}{{j-1}\choose i}\left(\frac{q-1}{q}\right)^i
tr\Big((g_1)_2(g_4)_{i+1}\Big)\; .$$
Hence, $tr(g_1g_2L_4)$ and $tr(g_1g_2L_5)$ consist of reduced traces that we
already evaluated, and yield no new
information.
$tr(g_1g_2L_6)$ gives rise to the new doubly connected term $tr(g_1g_2\,
g_4g_5)$.
In general,
$tr(g_1g_2L_m)$ yields, in addition to known terms, the new term
$tr\Big((g_1)_2\, (g_4)_{m-4}\Big)$.

For $tr(g_1g_2g_3L_m)$ we obtain from eq. \ref{equation:LkLm}
$$tr(g_1g_2g_3L_m)=\left(\frac{q}{q-1}\right)^2\Big(tr(L_4L_m)
-3tr(L_3L_m)+3tr(L_2L_m)\Big)$$
On the other hand,
\begin{eqnarray}
tr(g_1g_2g_3L_m) &=&
\frac{1}{q^{m-1}}tr\Big((g_1)_3(g_1)^{(m)}\Big)+
\frac{1}{q^{m-2}}tr\Big((g_1)_3(g_2)^{(m-1)}\Big)
\nonumber \\
&+&\frac{1}{q^{m-3}}tr\Big((g_1)_3(g_3)^{(m-2)}\Big)+
\frac{1}{q^{m-4}}tr\Big((g_1)_3(g_4)^{(m-3)}\Big)\nonumber \\ &+&
\sum_{i=5}^m\frac{1}{q^{m-i}}tr\Big((g_1)_3(g_i)^{(m-i)}\Big)
\nonumber
\end{eqnarray}
New reduced terms only appear for $m\geq 8$, consisting of a sequence of length
three ({\it i.e.,} a cycle of length {\it four}), and a disjoint sequence of
length $m-5$.
For $m=6,7$ we obtain reduced traces with two disjoint sequences of lengths
$3+1$ and
$3+2$, respectively, that have already been evaluated above.
The first four terms are
evaluated using eqs. \ref{equation:l3l}, \ref{equation:g1gi},
\ref{equation:ll},
and \ref{equation:l1l}, respectively.

 From eq. \ref{equation:g1g3}
\begin{equation}
\label{equation:g1g3L6}
tr(g_1g_3L_6)=tr(L_2L_4L_6)-\frac{q^2+1}{q-1}tr(L_4L_6)+
\frac{(q+1)^2}{q-1}tr(L_3L_6)-\frac{2q}{q-1}tr(L_2L_6)\; .
\end{equation}
Now,
\begin{eqnarray}
\label{equation:eighteen}
tr(g_1g_3L_6) &=& tr(g_1g_3g_5)+\frac{1}{q}tr\Big(g_1g_3(g_4)^{(2)}\Big)+
\frac{1}{q^2}tr\Big(g_1g_3(g_3)^{(3)}\Big)\nonumber \\ &+&
\frac{1}{q^3}tr\Big(g_1g_3(g_2)^{(4)}\Big) +
\frac{1}{q^4}tr\Big(g_1g_3(g_1)^{(5)}\Big)
\end{eqnarray}
Only the last two terms on the right hand side do not immediately follow from
identities previously derived.
Presenting them in a slightly generalized form we obtain
$$tr\Big(g_1g_3(g_2)^{(\ell)}\Big)=(q-1)(q^2+1)tr\Big((g_1)_3(g_4)^{(\ell-2)}\Big)+
(q-1)^2q\, tr\Big((g_1)_2(g_3)^{(2)}\Big)+
q^2tr\Big(g_1g_3(g_4)^{(2)}\Big)$$
and
\begin{eqnarray}
tr\Big(g_1g_3(g_1)^{(\ell)}\Big) &=&
(q-1)(q^2+1)(q^2-q+1)tr\Big((g_1)_3(g_4)^{(\ell-3)}\Big) \nonumber \\ &+&
(q-1)^2q(2q^2-q+2)tr\Big((g_2)_2(g_4)^{(\ell-3)}\Big)+
q^2(q^2-q+1)tr\Big(g_1g_3(g_4)^{(\ell-3)}\Big) \nonumber \\ &+&
q^2(q-1)^3tr\Big(g_3(g_4)^{\ell-3)}\Big)+
(q-1)q^3tr\Big(g_2(g_4)^{(\ell-3)}\Big) \nonumber
\end{eqnarray}
Hence, eq. \ref{equation:eighteen} provides the triply connected trace
$tr(g_1g_3g_5)$.

Similarly, replacing $L_6$ by $L_7$ in eq. \ref{equation:g1g3L6} we can derive
an expression for the reduced trace $tr(g_1 \, g_3\, g_5g_6)$.

The explicit expressions for all the reduced traces appearing up to $H_7(q)$,
in terms of traces of products of Murphy operators, are presented in Table 1.

\section{Evaluation of the traces of products of Murphy operators within
irreps of $H_n(q)$}

\hspace*{6 mm}
In the previous section we have  expressed the reduced traces
in terms of traces of products of Murphy operators.
We now present a recursive procedure that can be implemented to obtain the
latter traces.

For an irrep $\Gamma_n$ of $H_n(q)$ consider the set of irreps
$\Gamma_{n-1}\subset \Gamma_n$ of $H_{n-1}(q)$, obtained by eliminating
one box from $\Gamma_n$ in all possible ways.
Clearly,
$$tr(L_i)_{\Gamma_n}=\sum_{\Gamma_{n-1}\subset\Gamma_n}
tr(L_i)_{\Gamma_{n-1}}\;\;\;
i=2,\, 3,\, \cdots,\, n-1.$$
and
$$|\Gamma_n|=\sum_{\Gamma_{n-1}\subset\Gamma_n} |\Gamma_{n-1}|$$
where $|\Gamma_n|$ is the dimensionality of the irrep $\Gamma_n.$

\noindent
The trace of $L_n$ can now be evaluated very conveniently using
$$tr(L_n)_{\Gamma_n}=
\sum_{\Gamma_{n-1}\subset\Gamma_n} |\Gamma_{n-1}|\,
\{ \Gamma_n\setminus\Gamma_{n-1} \}_q\; .$$

We proceed to obtain the traces of
products of distinct, non-consecutive Murphy operators.
For
$\prod_{i=1}^{\ell} L_{\alpha_i}$,
$\alpha_1\geq 2$, $\alpha_{i+1}\geq\alpha_i+2 \, , \; i=1,\, 2,\, \cdots,\,
\ell-1$,
and $\alpha_{\ell}\leq n$,
we consider all the sequences of irreps leading to the desired irrep $\Gamma_n$
of
$H_n(q)$ from all possible irreps of $H_{\alpha_1-1}(q)$.
Let
$\Gamma_{\alpha_1-1}\subset\Gamma_{\alpha_1}\subset\Gamma_{\alpha_1+1}\subset\cdots
\subset\Gamma_n$ be such a sequence and let $\{\cdots\subset\Gamma_n\}$ denote
the
complete set of relevant sequences.
Obviously,
$$tr\left(\prod_{i=1}^{\ell} L_{\alpha_i} \right)_{\Gamma_n}=
\sum_{\{\cdots\subset\Gamma_n\}}
|\Gamma_{\alpha_1-1}| \prod_{i=1}^{\ell}
\{ \Gamma_{\alpha_i}\setminus\Gamma_{\alpha_i-1} \}_q \; .$$
This expression can be rewritten in the following form:

\noindent
if $\alpha_{\ell}<n$:
$$tr\left(\prod_{i=1}^{\ell} L_{\alpha_i}\right)_{\Gamma_n}=
\sum_{\Gamma_{n-1}\subset\Gamma_n}
tr\left(\prod_{i=1}^{\ell} L_{\alpha_i}\right)_{\Gamma_{n-1}}$$

\noindent
if $\alpha_{\ell}=n$:
$$tr\left(\prod_{i=1}^{\ell} L_{\alpha_i}\right)_{\Gamma_n}=
\sum_{\Gamma_{n-1}\subset\Gamma_n}
\{\Gamma_n\setminus\Gamma_{n-1}\}_q
tr\left(\prod_{i=1}^{\ell-1} L_{\alpha_i}\right)_{\Gamma_{n-1}}$$

Hence, given the traces of all Murphy operators up to $L_{n-2}$ and of all
products of sets of non-consecutive Murphy operators, corresponding to the
irreps
of $H_{n-1}(q)$, those corresponding to the irreps of $H_n(q)$ are readily
obtained.

The procedure presently described is ideally suited to implementation using
symbolic
programming.
Such an implementation was used to evaluate the traces of all
relevant products of Murphy operators in $H_6(q)$ and $H_7(q)$, in terms of
which
the corresponding reduced traces were obtained.
The results for $H_6(q)$ agree with those
in ref. \cite{King}, and the results for $H_7(q)$ are presented in table 2.
For $q=1$ they reduce to the characters of the symmetric group $S_7$.
As a further check we note that the traces of reduced elements consisting of
$\ell$
generators, that correspond to conjugate Young diagrams, are related by the
transformation $q^{\ell-i}\mapsto (-1)^{\ell}q^i$, $i=0,\, 1,\, \cdots,\,
\ell.$

\section{Conclusions}

\hspace*{6 mm}
In the present article we have demonstrated that the trace of an arbitrary
element of the Hecke algebra $H_n(q)$ can be reduced in a systematic way into
a linear combination of reduced traces, consisting of disjoint sequences
of consecutive generators. The reduced traces depend only on their cycle
structures, {\it i.e,} on the set of lengths of their disjoint sequences,
and can be expressed in terms of linear combinations of traces of products
of Murphy operators. The latter traces can be evaluated using a straightforward
recursive procedure.

One favourable feature of this approach is
that the expression for a given type of
reduced trace in terms of traces of products of Murphy opeartors does not
explicitly depend on $n$. Consequently, the present procedure is
particularly attractive when specific types of traces, {\it e.g.,} traces with
some
restriction on the set of disjoint cycles, are required to high order in $n$.
Thus, the present scheme may be of interest towards the investigation of the
high
$n$ limit of restricted classes of reduced traces, a problem of some current
interest
\cite{Kerov}.

\vspace{4mm}

\noindent
{\bf Acknowledgements}

Helpful discussions with Alain Lascoux, Bernard Leclerc,  Jean-Yves Thibon
and G\'erard Duchamp are gratefully acknowledged.

\newpage

\noindent
{\bf APPENDIX : The trace reduction procedure.}

In this appendix we show how to obtain the reduced form of the trace

\begin{eqnarray}
\label{equation:A1}
&tr&\Big((g_{1}...g_{m_{1}-1})(g_{m_{1}+1}...g_{m_{2}-1})...
(g_{m_{j}+1}...g_{p})(g_{k}...g_{p+r}...g_{k})\Big) \nonumber\\
&\equiv& tr\Big((g_{1})_{m_{1}-1}(g_{m_{1}+1})_{m_{2}-m_{1}-1}...
(g_{m_{j}+1})_{p-j}(g_{k})^{p+r-k+1}\Big)\nonumber \\
\end{eqnarray}
where
$(g_{1})_{m_{1}-1}(g_{m_{1}+1})_{m_{2}-m_{1}-1}...
(g_{m_{j}+1})_{p-j}$ is a reduced element consisting of $j+1$
connected sequences, or, in other words, a sequence $g_1g_2\cdots g_p$ with $j$
cuts
at $m_{i}$, $i=1,...,j$, and $(g_k)^{(p+r-k+1)}$,
$1\leq k \leq p+1$, is a term contributed by the Murphy operator $L_{p+r+1}$.
It is shown in section 6 that the reduction of eq. \ref{equation:A1}
is sufficient in order to express all reduced traces in terms
of the traces of products of Murphy operators.

A central role is played by the coefficients $f_{p}$ defined
through
$$g_i^p=f_pg_i+qf_{p-1}$$
They satisfy the recurrence relation
\begin{equation}
\label{equation:e3}
f_{p}=(q-1)f_{p-1}+qf_{p-2}
\end{equation}
with $f_{0}\equiv 0$ and $f_{1}=1$, that yields
\begin{equation}
\label{equation:e4}
f_{p}=\frac{q^p-(-1)^p}{q+1}= q^{p-1}-q^{p-2}+...+(-1)^{p-1}
\end{equation}
In terms of these we will obtain {\it f-expansions} leading to completely
reduced forms.

Systematic use will be made of the lemmas
\begin{equation}
\label{equation:e5}
g_{i}(g_{j}g_{j+1}...g_{k})=(g_{j}g_{j+1}...g_{k})g_{i-1}\;\;\;\;\;
(j < i \leq k)
\end{equation}
and
\begin{equation}
\label{equation:e6}
[g_{i},(g_{j}...g_{j+r}...g_{j})]=0\;\;\;\;\;\;(j<i<j+r)
\end{equation}

We start with the simplest form of eq. \ref{equation:A1}, allowing no cuts
($j=0$),
proceeding step-wise up in the order of
complexity.

Define
$$V_k^{(p,r)}\equiv tr\Big((g_{1}...g_{k-1})(g_{k+1}...g_{p})(g_{p+1}...
g_{p+r}...g_{p+1})\Big)
\;\; ;\;\;\;\;0 \leq k \leq p-1$$
where $(g_i)_0\equiv 1$ so that
$$V_0^{(p,r)}\equiv tr\Big((g_{1}...g_{p})(g_{p+1}...g_{p+r}...g_{p+1})\Big)$$
Using \ref{equation:e5} and \ref{equation:e6} one easily obtains the reduced
form
$$V_{k}^{(p,r)}=\sum_{\ell=0}^{r-1}{{r-1}\choose l}q^{l}(q-1)^{r-l-1}
tr\Big((g_{1}...g_{k-1})(g_{k+1}...g_{p+r-l})\Big)$$
For fixed $p$ and $r$, when there is no risk of confusion, we will write
$V_{k}$ for $V_{k}^{(p,r)}$. Consider now the case with $p-l+1$ overlapping
generators,
but no cut, namely
$$tr\Big((g_1...g_p)(g_l...g_{p+r}...g_l)\Big) \;\;\;\;(\ell\leq p)$$
and note that for $l=p+1$ one has just $V_{0}^{(p,r)}$, and for $l>p+1$ the
reduced
form is
obtained trivially.

For $p\geq 2l-1$, repeatedly using the lemmas \ref{equation:e5} and
\ref{equation:e6} along with cyclic
permutations within the trace, one obtains
\begin{eqnarray}
\label{equation:e10}
&tr&\Big((g_{1}...g_{p})(g_{l}...g_{p+r}...g_{l})\Big) \nonumber\\
&=&tr\Big((g_{p-2k}...g_{l-k})(g_{1}...g_{p})(g_{l-k}...g_{p})
(g_{p+1}...g_{p+r}...g_{p+1})\Big)\;\;\;\;\;\;\;\;\;(k\leq l-1) \nonumber\\
&=&tr\Big((g_{p-l+1}...g_{1})(g_{1}...g_{p-l+1})(g_{1}...g_{p})
(g_{p+1}...g_{p+r}...g_{p+1})\Big)
\end{eqnarray}
The final form is well defined for all $p\geq l$, but the
intermediate steps make the inequality $p\geq 2l-1$ necessary.

Using eq. \ref{equation:Lpp}
$$(g_{p-l+1}...g_{1})(g_{1}...g_{p-l+1})
=(q-1)\sum_{k=1}^{p-l+1}q^{k-1}(g_{k}...g_{p-l+1}...g_{k})+q^{p-l+1}$$
and, as follows from eq. \ref{equation:e5},
\begin{eqnarray*}
&tr&\Big((g_{k}...g_{p-l+1}...g_{k})(g_{1}...g_{p})
(g_{p+1}...g_{p+r}...g_{p+1})\Big) \\
&=&tr\Big((g_{1}...g_{p-l-k+2}...g_{1})(g_{1}...g_{p})(g_{p+1}...
g_{p+r}...g_{p+1})\Big)
\end{eqnarray*}
we obtain from \ref{equation:e10}
\begin{equation}
\label{equation:e14}
tr\Big((g_{1}...g_{p})(g_{l}...g_{p+r}...g_{l})\Big)
=(q-1)\sum_{r=1}^{p-l+1} q^{r-1}A_{p-l-r+2}^{(p,r)}+q^{p-l+1}V_{0}^{(p,r)}
\end{equation}
where
\begin{equation}
A_{s}^{(p,r)}\equiv tr\Big((g_{1}...g_{s}...g_{1})(g_{1}...g_{p})
(g_{p+1}...g_{p+r}...g_{p+1})\Big)
\end{equation}
Skipping details and suppressing the superscript $(p,r)$ we state
the recurrence relation
\begin{equation}
\label{equation:e15}
A_{s}=(q-1)sq^{s-1}V_{0}+q^{s}V_{s}
+(q-1)^{2}\sum_{r=1}^{s-1}rq^{r-1}A_{s-r}
\end{equation}
In the derivation an auxiliary construction
$$X_{(k,l)}\equiv tr\Big((g_{1}...g_{k})(g_{k+1}...g_{l}...g_{k+1})
(g_{k}...g_{p})(g_{p+1}...g_{p+r}...g_{p+1})\Big)\; ,$$
satisfying
$$X_{(k,l)}=(q-1)A_{l-k}+q X_{(k,l-1)} \; ,$$
was used.
Eq. \ref{equation:e15} is satisfied by
\begin{equation}
\label{equation:e18}
A_{s}=f_{2s}V_{0}+(q-1)\sum_{m=1}^{s-1}q^{m}f_{2(s-m)}V_{m}+q^{s}V_{s}
\end{equation}
This follows, by induction, noting that
\begin{eqnarray*}
A_{1}&=&tr\Big(g_{1}(g_{1}...g_{p})(g_{p+1}...g_{p+r}...g_{p+1})\Big)
\\
&=&(q-1)V_{0}+qV_{1}
\end{eqnarray*}
and that, using \ref{equation:e4} and summing the series in $k$, one obtains
the identity
\begin{equation}
\label{equation:Atwo}
(q-1)^2\sum_{k=1}^{p-1} k q^{k-1}f_{2(p-k)} + (q-1)pq^{p-1}=f_{2p} \; ,
\end{equation}
with the aid of which the induction is straightforward.
Finally, injecting \ref{equation:e18}
in \ref{equation:e14} and using the identity (obtained like
\ref{equation:Atwo}),
$$(q-1)\sum_{\ell=1}^p q^{\ell-1}f_{2(p+1-\ell)} + q^p = f_{2p+1}$$
one obtains,
\begin{equation}
\label{equation:e22}
tr\Big((g_{1}...g_{p})(g_{l}...g_{p+r}...g_{l})\Big)
= f_{2(p-l)+3}V_{0}+(q-1)\sum_{k=1}^{p-l+1} q^{k}f_{2(p-l-k)+3}V_{k}
\end{equation}
In the partial reduction \ref{equation:e10} we had to impose $p\geq 2l-1$.
For the complementary situation $(p< 2l-1)$ {\it the same final result},
eq. \ref{equation:e22},
can be obtained following a different route. This will now be briefly
indicated. One can show, using the lemmas along with cyclic permutations, that
\begin{eqnarray}
\label{equation:e23}
&tr&\Big((g_{1}...g_{p})(g_{l}...g_{p+r}...g_{l})\Big) \nonumber\\
&=& tr\Big((g_{1}...g_{2(p-l+1)})(g_{1}...g_{p})(g_{p+1}...g_{p+r}...
g_{p+1})\Big)
\end{eqnarray}
Define
$$B_{k}=tr\Big((g_{1}...g_{k})(g_{1}...g_{p})(g_{p+1}...g_{p+r}...g_{p+1})\Big)$$
The following recurrence relations have been obtained
\begin{equation}
\label{equation:e25}
B_{2(k+1)}=q^{k+1}V_{0}+(q-1)\sum_{r=0}^{k} q^{k-r} B_{2r+1}
\end{equation}
\begin{equation}
\label{equation:e26}
B_{2k+1}=q^{k+1}V_{k+1}+(q-1)\sum_{r=0}^{k} q^{k-r} B_{2r}
\end{equation}
In the derivation we used the auxiliary construction
$$Y_{(k,s)}\equiv tr\Big((g_{1}...g_{k-s})(g_{s+1}...g_{p})(g_{p+1}...
g_{p+r}...g_{p+1})\Big)\;\;\;\;\;(k\geq 2s+1, p>k)$$
that satisfies
$$Y_{(k,s)}=(q-1)B_{(k-2s-1)}+qY_{(k,s+1)}$$
and
$$B_{k}=(q-1)B_{k-1}+qY_{(k,1)} \; .$$

The solutions of \ref{equation:e25} and \ref{equation:e26} are
$$B_{2k}=f_{2k+1}V_{0}+(q-1)\Big(\sum_{m=1}^{k}q^{m}f_{2(k-m)+1}V_{m}\Big)$$
$$B_{2k+1}=f_{2k+2}V_{0}+(q-1)\Big(\sum_{m=1}^{k}q^{m}f_{2(k-m+1)}V_{m}\Big)+
q^{k+1}V_{k+1}$$

Since the right hand side of \ref{equation:e23} is just $B_{2(p-l+1)}$
we get eq. \ref{equation:e22} again,
as promised. Hence, the crucial data (both for $p\geq 2l-1$ and for $p<2l-1$)
determining the structure of the $f$-expansion, eq. \ref{equation:e22},
is the length of the
overlap $p-l+1$. This is reminiscent of, but distinct from, the result that
the reduced traces depend
only on the lengths of the disjoint cycles.

We now generalize to \ref{equation:A1} with a single cut, namely
$$tr\Big((g_{1}...g_{k-1})(g_{k+1}...g_{p})(g_{l}...g_{p+r}...g_{l})\Big)  $$

{\bf Case 1.} $l>k$

The generalization is trivial. The factor $(g_{1}...g_{k-1})$ commutes with
all the generators in the other two factors, that can be reduced essentially
as in the previous case (without cut), with $(g_{1}...g_{k-1})$ as a passive
spectator.

{\bf Case 2.} $l=k$

One has
\begin{eqnarray*}
&tr&\Big((g_{1}...g_{p})(g_{l}...g_{p+r}...g_{l})\Big) \\
&=&tr\Big((g_{1}...g_{l})(g_{l}...g_{p+r}...g_{l})(g_{l+1}...g_{p}
)\Big) \\
&=&(q-1)tr\Big((g_{1}...g_{l-1})(g_{l+1}...g_{p})(g_{l}...g_{p+r}
...g_{l})\Big)+q\, tr\Big((g_{1}...g_{p})(g_{l+1}...g_{p+r}...g_{l+1})\Big)
\end{eqnarray*}
Hence,
\begin{eqnarray}
\label{equation:e34}
&tr&\Big((g_{1}...g_{l-1})(g_{l+1}...g_{p})(g_{l}...g_{p+r}...g_{l})\Big)
\nonumber\\
&=&{1\over{q-1}}\Big\{tr\Big((g_{1}...g_{p})(g_{l}...g_{p+r}...g_{l})\Big)-
q\, tr\Big((g_{1}...g_{p})(g_{l+1}...g_{p+r}...g_{l+1})\Big\} \\
&=&f_{2(p-l-1)}V_{0}+(q-1)\sum_{k=1}^{p-l}q^{k}f_{2(p-l-k+1)}V_{k}+
q^{p-l+1}V_{p-l+1} \nonumber
\end{eqnarray}
where we used eqs. \ref{equation:e22} and \ref{equation:e3}.

\newpage

An evident generalization of eq. \ref{equation:e34} is
\begin{eqnarray}
\label{equation:e36}
&tr&\Big((g_{1}...g_{m-1})(g_{m+1}...g_{l-1})(g_{l+1}...g_{p})(g_{l}...
g_{p+r}...g_{l})\Big) \\
&=&{1\over(q-1)}\Big\{tr\Big((g_{1}...g_{m-1})(g_{m+1}...g_{p})(g_{l}...
g_{p+r}...g_{l})\Big) \nonumber\\
&&{\phantom{space is needed here}}
-q\, tr\Big((g_{1}...g_{m-1})(g_{m+1}...g_p)(g_{l+1}...g_{p+r}...
g_{l+1})\Big)\Big\} \nonumber
\end{eqnarray}
An analogous result holds quite generally for $k=m_{j}$ in \ref{equation:A1}.

The results \ref{equation:e34}, \ref{equation:e36} and their evident
successive generalizations will be
used repeatedly for the cases to follow.

{\bf Case 3.} $l<k$

Now using lemma \ref{equation:e6}
\begin{equation}
\label{equation:e37}
tr\Big((g_{1}...g_{k-1})(g_{k+1}...g_{p})(g_{l}...g_{p+r}...g_{l})\Big)=
tr\Big((g_{1}...g_{k-1})(g_{l}...g_{p+r}...g_{l})(g_{k+1}...g_{p})\Big)
\end{equation}
One can reduce the first two factors on the right hand side of eq.
\ref{equation:e37},
leading to \ref{equation:e22},
with $(g_{k+1}...g_{p})$ as spectator. This gives
\begin{eqnarray*}
&tr&\Big((g_{1}...g_{k-1})(g_{k+1}...g_{p})(g_{l}...g_{p+r}...g_{l})\Big)
\\
&=&f_{2(k-l)+1}U_{0}+(q-1)\sum_{s=1}^{k-l}q^{s}f_{2(k-l-s)+1} U_{s}
\end{eqnarray*}
where
\begin{eqnarray}
\label{equation:e38}
U_{m}&=&tr\Big((g_{1}...g_{m-1})(g_{m+1}...g_{k-1})(g_{k}...g_{p+r}...
g_{k})(g_{k+1}...g_{p})\Big) \nonumber\\
&=&tr\Big((g_{1}...g_{m-1})(g_{m+1}...g_{k-1})(g_{k+1}...g_{p})(g_{k}...
g_{p+r}...g_{k})\Big)
\end{eqnarray}
In particular
$$U_{0}=tr\Big((g_{1}...g_{k-1})(g_{k+1}...g_{p})(g_{k}...g_{p+r}...g_{k})\Big)$$
has already been studied in case 2.
\begin{eqnarray*}
U_{1}&=&tr\Big((g_{2}...g_{k-1})(g_{k}...g_{p+r}...
g_{k})(g_{k+1}...g_{p})\Big) \\
&=&tr\Big((g_{1}...g_{k-2})(g_{k}...g_{p-1})(g_{k-1}...
g_{p+r-1}...g_{k-1})\Big)
\end{eqnarray*}
is also directly reduced to the preceding case.

For $U_{m}(m>1)$ one first uses eq. \ref{equation:e36} to obtain traces with a
single cut. After this step, since $s\leq k-l<k$, all the cuts in the traces
$U_s$
in eq. \ref{equation:e38} lie below the overlap, as for case 1. Hence, one can
at that stage obtain complete reduction.

{\it This process indicates clearly the generalization of the iterative
process necessary for multiple cuts, as in \ref{equation:A1}}.
One starts by shifting the
factor $(g_{m_{j}+1}...g_{p})$ on the right and uses for the remaining
factors the reduction procedure for $j-1$ cuts with known results.
Finally, one proceeds as for the passage from $j=0$ to $j=1$.

For the simple case 2 the $f$'s were recombined to obtain the simple form
(35). A possible generalization of such combinatorics for arbitrary cuts is
beyond the scope of this paper.

To sum up, we have the complete machinery permitting a multiple stage
$f$-expansion of \ref{equation:A1} to its reduced form.
A whole set of interesting algebraic
structures were discovered in the process.

\end{document}